\documentclass[pdflatex,sn-mathphys]{sn-jnl}

\jyear{2021}%

\theoremstyle{thmstyleone}%
\usepackage{subcaption}
\usepackage{support-caption}
\usepackage{silence}
\usepackage{url}
\theoremstyle{thmstyletwo}%
\theoremstyle{thmstylethree}%
\usepackage{amsmath}
\begin{document}

\title[Monu Yadav and Laxminarayan Das]{Identification of storm eye from Satellite image data using fuzzy logic with machine learning}


\author[1]{\fnm{Monu} \sur{Yadav}}\email{yadavm012@gmail.com}

\author*[1]{\fnm{Laxminarayan} \sur{Das}}\email{lndas@dce.ac.in}


\affil[1]{\orgdiv{Department of Applied Mathematics}, \orgname{Delhi Technological University}, \orgaddress{ \state{New Delhi}, \country{India}}}




\abstract{This research presents a study of a unique technique for identifying storm eye that is based on fuzzy logic and image processing with the help of cloud images. Fuzzy logic is a term that refers to complicated systems with unclear behaviour caused by a number of different circumstances. It provides the ability to model the dynamic behavior of the storm and determines the location of the best eye in an area of interest. After that, image processing is applied to enable accurate eye positioning based on the search results.  The experimental results are analyzing the storm eye position with approxiamtely $98\%$ accurate compared to the India meteorological department provided best track data and Cooperative Institute for Meteorological Satellite Studies provided Advances Dvorak Technique data. As a result, the identification of storm's eye location using this technique can be found to improve significantly. Using the present technique, it is possible to determine the eye entirely automatically, which replacing the manual method that has been employed in the past.}




\maketitle

\section{Introduction}\label{sec1}
Along with gales, rainstorm, and storm surges, tropical cyclones (TCs) are considered extreme weather events that have the potential to cause massive damage for coastal areas. Meteorologists and warning centers have been studying the atmosphere, the atmospheric boundary layer, and the air-sea interface since the $19^{th}$ century and have made significant advances in observational technology, intensification physics, the ocean responses, and the interaction of these areas. In particular, predictions of TCs and their genesis, intensities, and risks remain quite difficult because of climate and geographical variation as well as  the climatological characteristics of oceans. As a result, a single forecasting technique will not be able to give results adequately and accurately. Additionally, there are multiple key characteristics of a storm, including wind speed, wind direction \cite{bib3,bib4,bib5,bib6,bib7,bib8,bib9,bib10,bib11}, pressure system \cite{bib3,bib4,bib7,bib12,bib13}, temperature \cite{bib4,bib5,bib11,bib12,bib14,bib15,bib16,bib27}, moisture precipitation \cite{bib4,bib12}, humidity\cite{bib12} and cloud shape \cite{bib3,bib45,bib17,bib18,bib14,bib15,bib19,bib20,bib21}. All of these factors affect prediction accuracy and can make forecasting more difficult.

Tracking the characteristics of prior information has led to the invention of several types of storm forecasting. Harr et al., 1995 \cite{bib28} said that the tropical cyclone (TC) variability is defined by the wind circulation patterns within tropical clusters and based on the storm trajectory tracking. Fang et al., 2015 \cite{bib29} presented a time-series approximation technique  for the prediction of a TC in the Western North Pacific was used for the development of the modified A-LTK (Approximation with use of Local features at Thinnedout Keypoints) technique for a multidimensional time-series similarity calculation method. In this technique, the time-series data is reduced to a smaller number of time points and the nearby time points are used to create a feature vector at the thinned-out time point. Using wind data provided by Joint Typhoon Warning Center and National Hurricane Center, Knaff et al., 2009 \cite{bib30} analyzed how peak intensity and lifetime maximum intensity could affect the recurrence of TCs. In several studies, eye location has been obtained from the ocean surface wind vectors extracted from SeaWind instrument from the satellite \cite{bib5,bib6,bib9,bib31,bib17,bib18}. Based on wind data, Artificial Intelligence Technique's Support Vector Machines were used to identify cyclones fully automatically \cite{bib6}, and Fuzzy C-Means clustering was used to identify cyclones objectively \cite{bib31}.

As an alternative source of information, cloud images can be used to determine the patterns of storm \cite{bib3,bib20,bib22,bib23,bib24,bib25,bib27} or the location of the storm eye\cite{bib17,bib18,bib14,bib15,bib19,bib21,bib26}. A neural network-based model and elastic pattern matching from the predefined TC technique known as Dvorak analysis were used for automatic TC pattern segmentation in \cite{bib3,bib16}, while Lee et al., 2000 \cite{bib27} proposed a similar idea using the Neural Oscillatory Elastic Graph Matching Model (NOEGM) for TC pattern identification. The Dvorak-based TC pattern matching algorithms underwent further development during \cite{bib22, bib23, bib24}. The algorithm determined the axisymmetry of the shape and flexible cloud cluster in TCs. Additionally, the algorithm predicted wind speed, storm trajectories and the expansion of TCs by the gradient of the shape and the Deviation-Angle Variance (DAV) of satellite infrared. Warunsin et al., 2015 \cite{bib19} proposed an automated approach to detect the center of a TC by fitting logarithmic spiral patterns to the image enhancement and the center of spiral has been assigned to the TC center. Although, this method may not be accurate in areas with ambiguous boundaries. As a remedy to this problem, mountain-climbing search was proposed by Bai et al., 2011 \cite{bib14}. The method was nevertheless technically challenging because it needed a lot of searching time and occasionally ran into problems with the conclusion. A proper termination condition for the searching algorithm should be taken into account to fix this issue. Due to this, Bai et al., 2012 \cite{bib15} introduced important requirements for the search algorithm's termination condition. An ant colony optimization approach was employed to advance the development of storm locating, and the Spiral Curve Model (SCM) was used to learn from the storm feature. However, it is still challenging to anticipate the exact storm location, particularly when there are unclear or many storm eyes.

The use of these techniques has been satisfying the needs of the Meteorological community for quiet a few decades, but now that there is a large volume of satellite observational, some advance mathematical agents as well as image processing, it seems that there is still room for improvement. As a mathematical agent, fuzzy is used for forecasting TC. The concept of fuzzy-based decision making allows for expert support at a linguistic level while evaluating data to learn statistical characteristics or to improve the rule. To improve detection results, it is necessary to have domain knowledge of cloud images of the storm. TC forecasting was proposed as a case-based fuzzy multi-criteria decision support system in \cite{bib32}. In case-based reasoning, we collected historical observation data for a ten-year period to develop fuzzy multi-criteria decision support. The Fuzzy Inference System (FIS) techique provide an advantange of reducing large volumes of data into a smaller number of the best cases, which then allows for further filtering and evaluation of the results for use by a forecaster in real time. Nevertheless, the accuracy of the algorithm varied since location of the storm center in a geographical region still played an important role in the forecast. Therefore, in order to handle the uncertainties of the storms, we have introduced a fuzzy-based storm eye identification system designed specifically for improving tracking of storm centers. An image processing technique for predicting TC characteristics can be used in a variety of ways. TC position and intensity are estimated reasonably accurately by passive and active microwave sensors such as the TRMM Microwave Imager and wind scatterometer \cite{bib36,bib37}. As a result of the high temporal sampling and wide coverage provided by the sensors, geostationary satellite IR images provide the best information for locating and estimating the intensity of TCs in real-time \cite{bib38}. TC analysts can now make use of dynamic or quickly evolving features in TC because geostationary satellites usually provide half hourly observations. As with intensity estimation, positioning procedures are based on the results of TC imagery conducted years ago \cite{bib39,bib40}. Some of the most important techniques for tracking TCs includes wind field analysis and pattern matching \cite{bib41,bib42}, the TREC (tracking radar echoes by correlation) algorithm \cite{bib43} or automatic cloud feature tracking \cite{bib44}. Several researchers, using radar data and IR images, have attempted to determine the eye of TC in recent years. A method to determine the approximate eye of any TC was developed by Liu et al. (2006) \cite{bib35} using the centre of the biggest circle contour within the TC in an infrared image. Yurchak et al., 2007 \cite{bib34} demonstrated that TC spiral bands can be represented analytically by a combination of hyperbolic and logarithmic spirals. Yan et al., 2005 \cite{bib33} and Wong et al., 2007 \cite{bib45} determined the eye of TC using motion field analysis. Compared with earlier eye determined studies, this one presents a new technique that is more unique. An automated method is used in this paper to determine the cyclone's eye. Fuzzy logic and Image processing are included in this. As part of fuzzy logic, we construct 16 If-Then rules of Multiple Input with Single Output system, define Membership function for Moisture Density, Wind speed, Pressure Drop at the Center, Estimated central pressure and defuzzification with Center of Area method to obtain crisp output. As part of Image processing, the techniques used include removal of noise from images by morphological operations as well as determining of center of gravity/mass of the TC image. Compared to other method, the proposed method shows improved performance. What are the important of fuzzy logic in image processing? Fuzzy logic used to manage the uncertainty and imperfection of an image. 
A description of data analysis is given in Section \ref{sec2}. The proposed method is presented in Section \ref{sec3}. Experimental results are presented in Section \ref{sec4} and compared with other method. Finally summary and conclusion are given in Section \ref{sec5}.

\section{Analyzing the Data} \label{sec2}
Typically, the eye of the storm has an important factor in TC tracking. Six-hourly IR images were analysed in this study. Satellite images generally will show TCs as clusters of well-organized clouds, although clouds clusters will lose structure during the formation process, causing the system to be difficulty to identify as a TC. The storm eye is difficult to locate due to incomplete and insufficient data.

There are IR images recorded from tropical storms and typhoon occurring in the North Indian Ocean and Bay of Bengal. For weather forecasting, IR images are widely used. It is generally accepted that several agencies collect the best weather tracks of a storm, such as India Meteorological Department (IMD), Joint Typoon Warning Center (JTWC), and Cooperative Institute of Meteorological Satellite Studies (CIMSS). However, the data for the best track recordings from each organization differ quite a bit from one another. Performance estimation will be affected by these track errors.
\section{A Possible Location of the Storm Eye} \label{sec3}
\begin{figure}
    \centering
    \includegraphics[width = 1\textwidth]{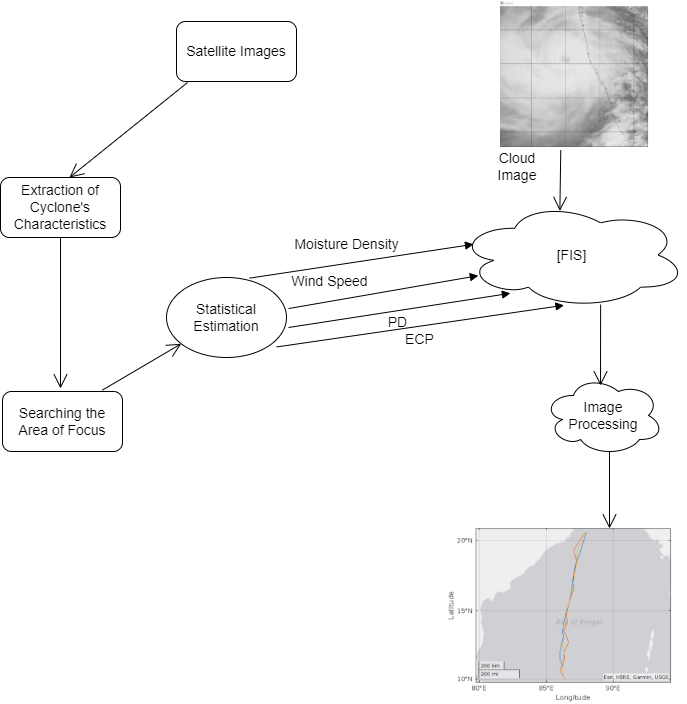}
    \caption{Proposed System. PD stand for Pressure drop at the center; ECP stand for Estimate central pressure.}
    \label{1}
\end{figure}

To identify a storm, it is usually necessary to extract its characteristics (moisture density, wind speed, pressure drop at center, and estimated central pressure). From satellite data it is possible to see the storm's eye, which is one of the most widely observed storm features. In recent years, researcher have found many ways to detect the storm's eye, using artificial intelligence technique such as support vector machine, feature vector, genetic algorithms, fuzzy C-means, etc. The purpose of this study is to determine the storm eye using a novel FIS definition and then to further refine the location of the storm eye using image processing. The procedure is shown in figure \ref{1}.
\subsection{Searching the area of focus}
\begin{figure}[h]
	\centering
	\includegraphics[width=1\textwidth]{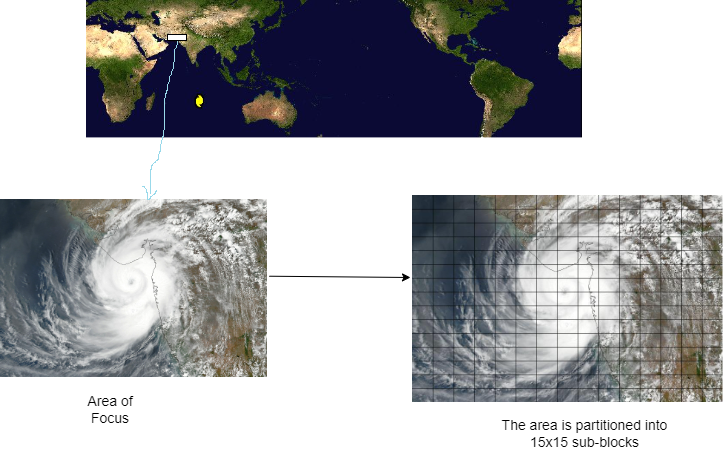}
	\caption{The process used to get the area of focus}
	\label{123}
\end{figure}

It was suggested that the storm's starting wind speed woud be 17 m/s. However, the slower wind speed must be taken into consideration when considering the gestation duration. The center of area of focus was used to determine the initial of the storm's eye. The storm's approximate diameter ranges from 100 to 2000 km. The area of focus is $13.2 \times 13.2$ degree ($1400 \times 1400 $ $km^{2}$) from the storm's initial center in order to cover the whole storm area. The area of focus is partitioned into 225 ($15 \times 15$) sub-blocks as shown in figure \ref{123}.
\subsection{ Fuzzy Interference System (FIS)}
\subsubsection{Fuzzification}
\begin{figure}
    \centering
    
    \begin{subfigure}[b]{0.49\textwidth}
         \centering
         \includegraphics[width=1\textwidth]{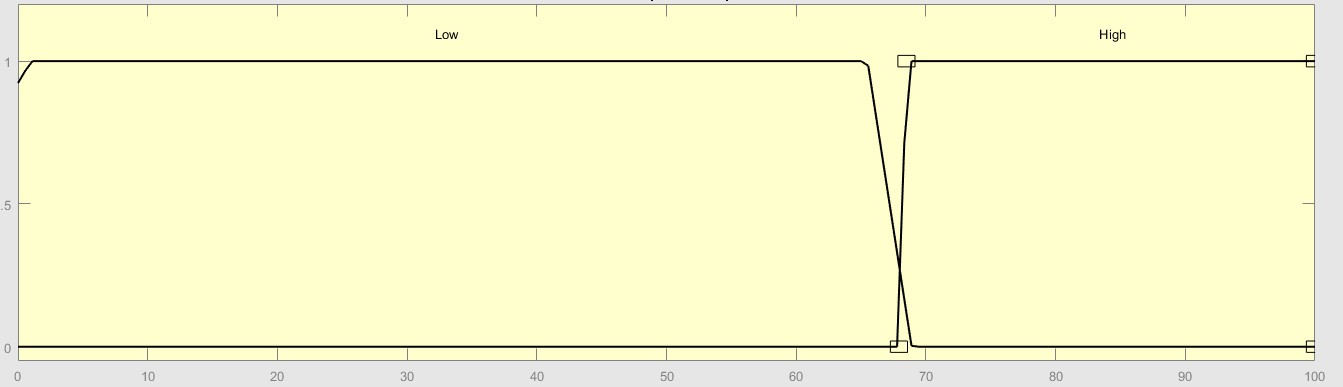}
         \caption{ Moisture Density}
         
     \end{subfigure}
     \hfill
     \begin{subfigure}[b]{0.49\textwidth}
         \centering
         \includegraphics[width=1\textwidth]{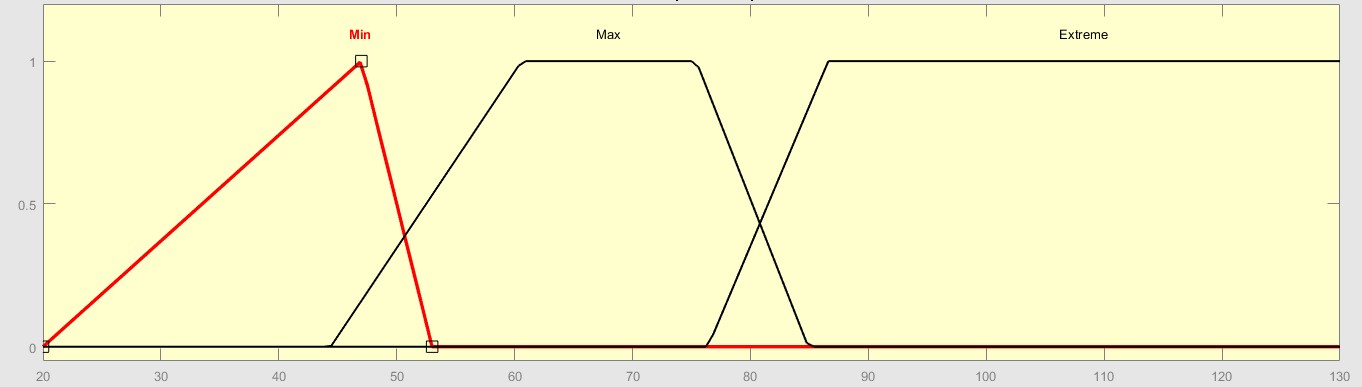}
         \caption{Wind Speed}
       
     \end{subfigure}
     \hfill
     \smallskip
     \begin{subfigure}[b]{0.49\textwidth}
         \centering
         \includegraphics[width=1\textwidth]{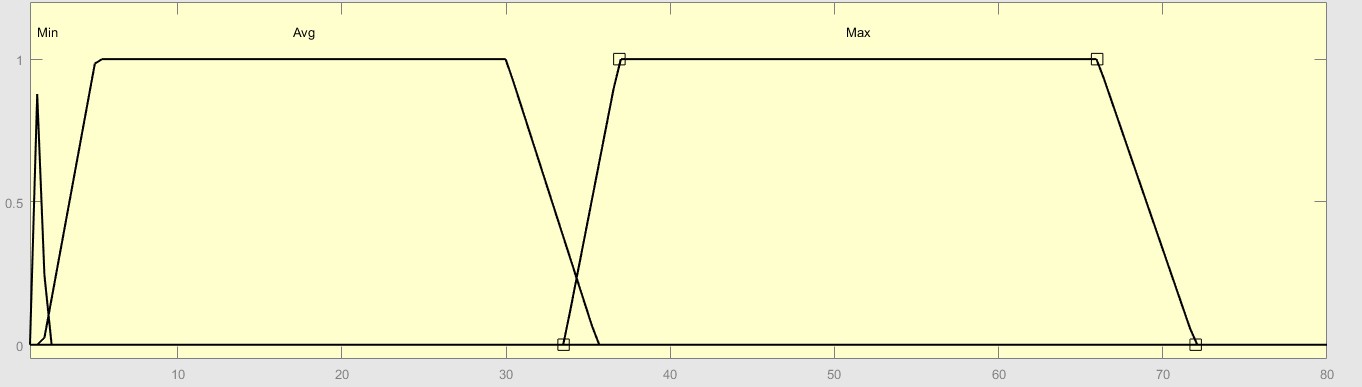}
         \caption{Pressure Drop at the Center}
         
     \end{subfigure}
     \hfill
     \begin{subfigure}[b]{0.49\textwidth}
         \centering
         \includegraphics[width=1\textwidth]{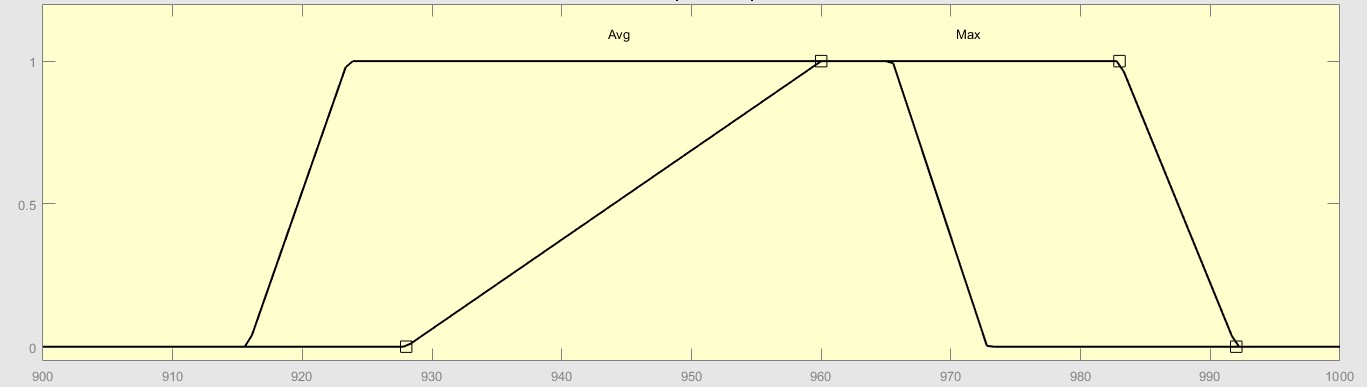}
         \caption{Estimated Central Pressure}
         
     \end{subfigure}
     \hfill
     \caption{Input provide to the proposed system are fuzzy in nature}
     \label{2}
\end{figure}

In this paper, a standard mamdani FIS \cite{bib1} was proposed to locate storms with expert knowledge. Based on the inputs such as moisture density, wind speed, pressure drop at the center, and estimated central pressure, Membership Function (MF) are constructed (Figure \ref{2}). A detailed description of meteorological input can be found in Table \ref{tab1}. The result is determined for each 6-hourly observation (03, 09, 18, 24 UTC) throughout the tropical cyclone's life period (we denoted the life period of tropical cyclone when its intensity is between 20 to 130 kts).

\begin{table}[ht]
    \centering
    \caption{A description to the inputs to the proposed FIS method}
    \begin{tabular}{|c|c|c|c|c|}
\hline
     Name & Symbol & Range &Unit & Membership Function$\footnotemark[2]$  \\
     \hline
    Moisture Density & D & 1-100 &\% &Low, High \\
    \hline
    Wind Speed & W &20-130 &Knots(kts)&Min, Max, Extreme \\
    \hline
    Pressure Drop at the Center&PD&  1-80& hPa$\footnotemark[1]$ & Min, Avg, Max \\
    \hline
    Estimated Central Pressure & EP & 900-1000 & hPa & Avg, Max \\
    \hline
    \end{tabular}
    \footnotetext[1]{Hectopascal is denoted as hPa, which is the SI unit of pressure} 
    \footnotetext[2]{Low, High, Min., Max., Avg., extreme are not a value but an interval.}
   
    \label{tab1}
\end{table}

\subsubsection{Inference}
\begin{figure}
    \centering
    \includegraphics[width = 1\textwidth]{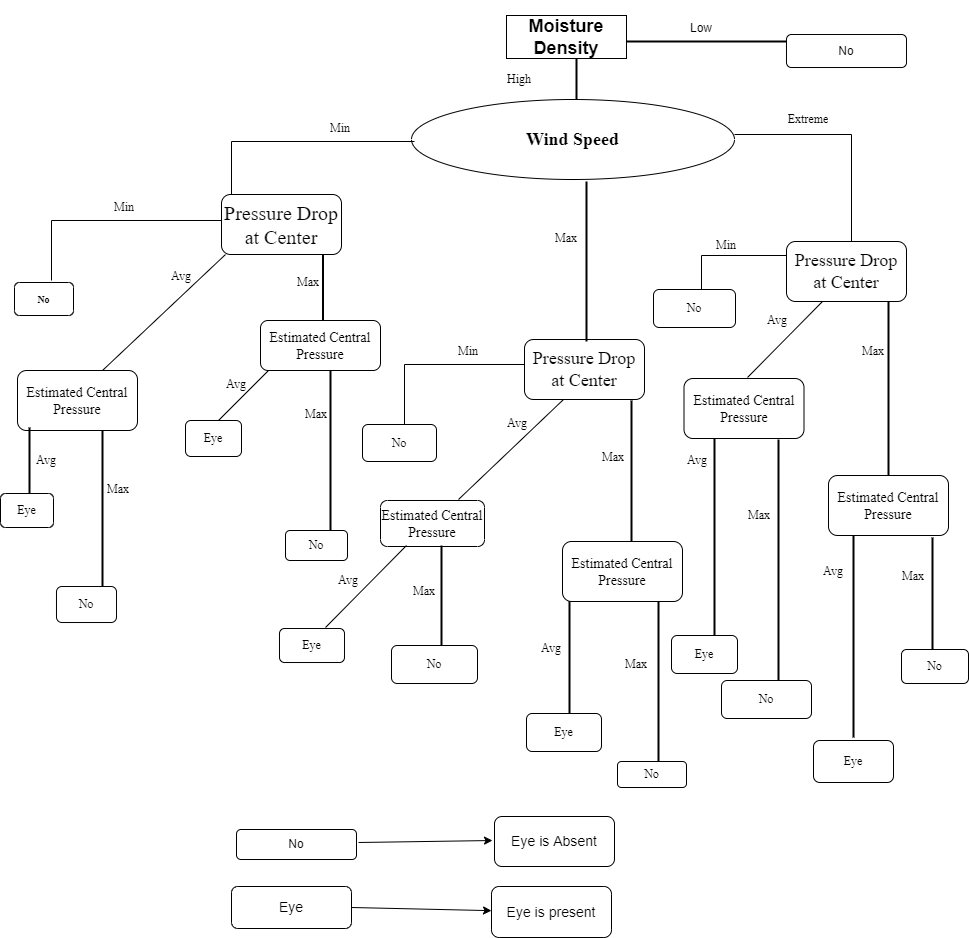}
    \caption{Decision tree for the proposed method}
    \label{3}
\end{figure}

Fuzzy inference rules are generated for determining the most likely location of the storm eye. In Figure \ref{3} , the proposed FIS system is shown as a decision tree. We construct 16 fuzzy IF-THEN based on Multi-Input-Single-Output systems are as follows:

\begin{enumerate}
    \item If (MoistureDensity is low) then (Eye is absent).
    \item  If (MoistureDensity is High) and (WindSpeed is Min) and (PressureDrop is Min) then (Eye is absent)                                          
    \item  If (MoistureDensity is High) and (WindSpeed is Min) and (PressureDrop is Avg) and (EstimatedCentralPressure is Avg) then (Eye is absent)      
    \item  If (MoistureDensity is High) and (WindSpeed is Min) and (PressureDrop is Avg) and (EstimatedCentralPressure is Max) then (Eye is present)        
    \item  If (MoistureDensity is High) and (WindSpeed is Min) and (PressureDrop is Max) and (EstimatedCentralPressure is Avg) then (Eye is present)        
    \item  If (MoistureDensity is High) and (WindSpeed is Min) and (PressureDrop is Max) and (EstimatedCentralPressure is Max) then (Eye is absent)      
    \item If (MoistureDensity is High) and (WindSpeed is Max) and (PressureDrop is Min) then (Eye is ansent)                                         
    \item  If (MoistureDensity is High) and (WindSpeed is Max) and (PressureDrop is Avg) and (EstimatedCentralPressure is Avg) then (Eye is present)        
    \item  If (MoistureDensity is High) and (WindSpeed is Max) and (PressureDrop is Avg) and (EstimatedCentralPressure is Max) then (Eye is present)        
    \item  If (MoistureDensity is High) and (WindSpeed is Max) and (PressureDrop is Max) and (EstimatedCentralPressure is Avg) then (Eye is present)       
    \item  If (MoistureDensity is High) and (WindSpeed is Max) and (PressureDrop is Max) and (EstimatedCentralPressure is Max) then (Eye is absent)    
    \item  If (MoistureDensity is High) and (WindSpeed is Extreme) and (PressureDrop is Min) then (Eye is absent)                                       
    \item  If (MoistureDensity is High) and (WindSpeed is Extreme) and (PressureDrop is Avg) and (EstimatedCentralPressure is Avg) then (Eye is present)   
    \item  If (MoistureDensity is High) and (WindSpeed is Extreme) and (PressureDrop is Avg) and (EstimatedCentralPressure is Max) then (Eye is absent) 
    \item  If (MoistureDensity is High) and (WindSpeed is Extreme) and (PressureDrop is Max) and (EstimatedCentralPressure is Avg) then (Eye is present)   
    \item  If (MoistureDensity is High) and (WindSpeed is Extreme) and (PressureDrop is Max) and (EstimatedCentralPressure is Max) then (Eye is absent) 
\end{enumerate}
\begin{figure}
    \centering
    \includegraphics[width = 1\textwidth]{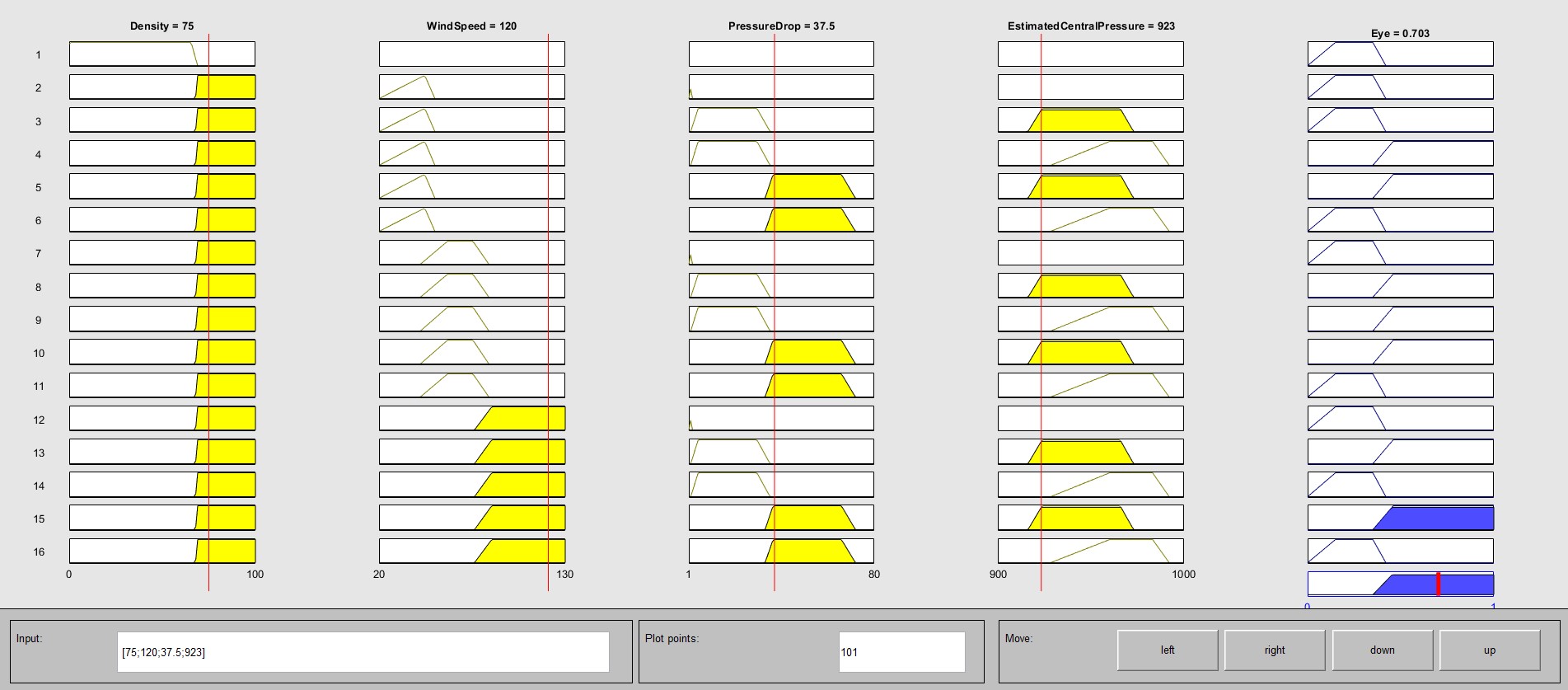}
    \caption{The visualisation of storm eye location using FIS taking an example. Unit of all parameter is mention in table \ref{tab1}. Numbering 1 to 16 are the if-then rule number. }
    \label{4}
\end{figure}
\subsubsection{Defuzzification}
The Center of Area method is the most general form of defuzzification. However, there are several other methods such as the Smallest of Maximas, the Largest of Maximas, and then Mean of Maximas. In Figure \ref{4}, we illustrate how FIS can help locate storm eyes is shown as taking an example. Moisture density is 75\%, wind speed is 120, pressure drop at the center is 37.5, and the estimated central pressure is 923 with the crisp inputs of density, wind speed, pressure drop at the center, and the estimated central pressure. We performed COA defuzzification and have obtained 0.703 as our crisp output. It is then concluded that the storm center is indeed within the identified block and the location of eye may subsequently be determined. Image processing should be carried out for eye location refinement in order to achieve an accurate location of the eye, as explained next.

\subsection{Image processing use for refinement of location of the Eye }
In the Image processing process, we refine the location of Eye in order to gain a more accurate forecast, as we known that Eye is one of the main characteristic to predict cyclones. Once the area of focus has been achieved, all of the images have equal pixel sizes. Change these images from grayscale to binary because binary images can be analysed more easily as they are only black and white. After that, morphological techniques such as dilation and erosion are used to remove noise from the images.

In standard practice, the position of the eye of the storm acts as the center of the storm. We estimate the center of mass in the images to be the storm center.
 To calculate the Center of gravity ($c_x, c_y$) for the image array of $m \times n$ order as follows :
\begin{equation*}
	c_{x} =\frac{\sum_{m}^ {i=1} \sum_{n}^{j=1} i \times C_{ij}}{\sum_{m}^ {i=1} \sum_{n}^{j=1}  C_{ij}} 
\end{equation*}
\begin{equation*}
	c_{y} =\frac{\sum_{m}^ {i=1} \sum_{n}^{j=1} j \times C_{ij}}{\sum_{m}^ {i=1} \sum_{n}^{j=1}  C_{ij}}
\end{equation*}
Among other factors, m refers to the number of horizonal pixels, n to the
number of vertial pixels. $C_{ij}$ represents the ($i^{th}$, $j^{th}$), element of an image.


\section{Result} \label{sec4}
Two of the TC cases have been tested with the algorithm discussed above, namely Tauktae, and Yaas. The satellite-derived IR images of these cyclones have been processed and the FIS method employed for determined the eye of TC and Image processing for refinement.

\begin{figure}
    \centering
    
    \begin{subfigure}[b]{0.47\textwidth}
        \includegraphics[width = 1\textwidth]{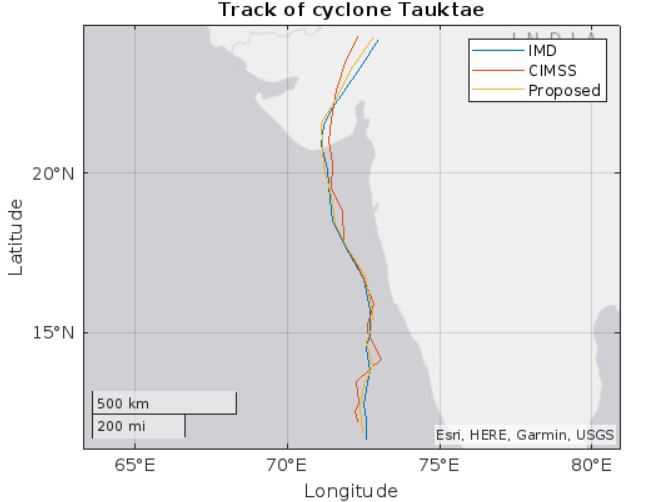}
    \caption{Track of cyclone Tauktae by IMD, CIMSS, Proposed Method}
    \end{subfigure}
    \hfill
    \begin{subfigure}[b]{0.47\textwidth}
        \includegraphics[width = 1\textwidth]{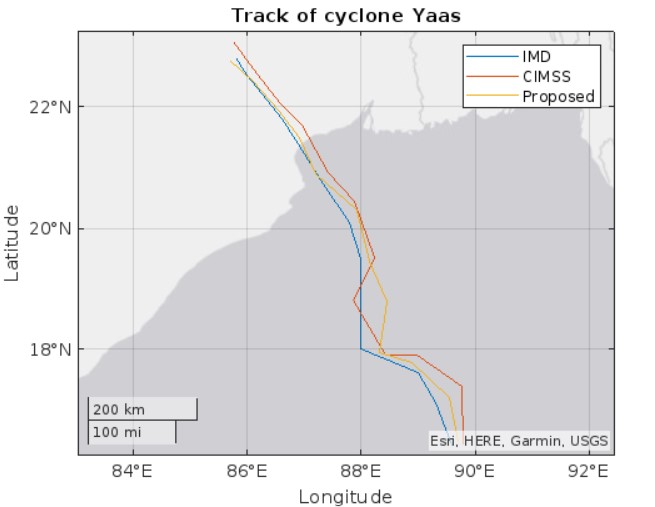}
    \caption{Track of cyclone Yaas by IMD, CIMSS, Proposed Method}
    \end{subfigure}
    \caption{Track of Two Tropical cyclones using proposed method and compare with IMD and CIMSS}
    \label{7}
\end{figure}
Figure \ref{7} illustrate the track of two TCs derived using the discussed method to determine the eye position, respectively. Comparing the cyclone track obtained from the discussed method with data from CIMSS and IMD.

Why do we compare TCs' tracks to validate our method? The TC's eye's motion is its track. The general circulation of the Earth's atmosphere has a significant role in their motion.

The error percentage of TC Tauktae between IMD, CIMSS, and the proposed method is shown in Table \ref{tab2}. We find that the difference between IMD and CIMSS is approximately 0.3 to 5 degree in latitude and 0.1 to 1 degree in longitude. In contrast, the difference between IMD and the proposed method is approximately 0.1 to 2 degree in latitude and less than 0.4 degree in longitude.

Table \ref{tab3} shows the TC Yaas error percentage between the proposed method, IMD, and CIMSS. We find that the IMD and CIMSS datasets differ by roughly 0 to 2 degrees in latitude and less than 1 degree in longitude. Conversely, the proposed method and IMD differ by about 0 to 1 degree in latitude and less than 0.5 degree in longitude.

An algorithm's accuracy is influenced by several factors, including the focus area, noise removal algorithm, and image processing algorithm. It is important that the center of TC should lie within the focus area in order for the algorithm to produce accurate results.

\begin{table}
	\centering
	\begin{tabular}{c c | c c| c c| c c c c }
		\hline
		& & & & & & & Error& Percentage& \\
		IMD & & CIMSS & & Proposed & &Lat. & Long. & Lat. & Long. \\
		Lat. & Long. & Lat. & Long. & Lat. & Long & A&A &B &B \\
		\hline
		11.6&72.&12.16&72.29&11.8&72.48&4.82&-0.427 &1.72&-0.16\\
		12.2&72.6&12.52&72.24&12.25&72.42&2.62&-0.49 &0.4&-0.24\\
		12.7&72.5&12.75&72.35&12.65&72.4&0.39&-0.2&-0.39&-0.13\\
		13.2&72.6&13.42&72.26&13.3&72.5&1.66&-0.46&0.75&-0.13\\
		13.8&72.7&14.16&73.1&13.95&72.8&2.6&0.55&1.08&0.13\\
		14.5&72.6&15.03&72.64&14.7&72.6&3.65&0.05&1.37&0\\
		15&72.7&15.25&72.65&15.1&72.75&1.66&-0.06&0.66&0.06\\
		15.7&72.7&15.9&72.83&15.67&72.78&1.27&0.17&-0.19&0.11\\
		16.7&72.5&16.9&72.43&16.8&72.55&1.19&-0.09&0.59&0.06\\
		17.5&72&17.87&71.86&17.7&71.95&2.11&-0.19&1.14&-0.06\\
		18.5&71.5&18.83&71.83&18.46&71.56&1.78&0.46&-0.21&0.08\\
		19.2&71.4&19.51&71.46&19.3&71.44&1.61&0.08&0.52&0.05\\
		20.1&71.3&20.3&71.5&20.15&71.2&0.99&0.28&0.24&-0.14\\
		20.9&71.1&20.98&71.37&20.85&71.15&0.38&0.37&-0.23&0.07\\
		21.5&71.2&21.65&71.46&21.55&71.1&0.69&0.36&0.23&-0.14\\
		22&71.5&22.54&71.61&22.25&71.55&2.45&.15&1.13&0.06\\
		23.1&72.3&23.48&71.95&23.15&72.05&1.64&-0.48&.21&-0.34\\
		24.1&73&24.2&72.36&24.15&72.85&0.41&-0.87&0.2&-0.21\\
		\hline

	\end{tabular}
\caption{Result of comparation with reference of TC, TAUKTAE. A stands for error percentage between IMD and CIMSS; B stands for error percentage between IMD and Proposed Method}
\label{tab2}
\end{table}

\begin{table}
	\centering
	\begin{tabular}{c c | c c| c c| c c c c }
		\hline
		& & & & & & & Error& Percentage& \\
		IMD & & CIMSS & & Proposed & &Lat. & Long. & Lat. & Long. \\
		Lat. & Long. & Lat. & Long. & Lat. & Long & A&A &B &B \\
		\hline
		16.4&89.6&16.41&89.78&16.45&89.7&0.06&0.2&0.30&0.11 \\
		17.1&89.3&17.37&89.76&17.19&89.55&1.57&0.51&0.52&0.27 \\
		17.6&89&17.89&88.97&17.78&88.9&1.64&-0.03&1.02&-0.11 \\
		18&88&17.9&88.43&17.95&88.32&-0.55&0.48&-0.27&0.36 \\
		18.7&88&18.82&87.87&18.79&88.45&0.64&-0.14&0.48&0.51 \\
		19.5&88&19.53&88.24&19.45&88.15&0.15&0.27&-0.25&0.17 \\
		20.1&87.8&20.44&87.88&20.3&87.92&1.69&0.09&0.99&0.13 \\
		20.8&87.3&20.93&87.42&20.89&87.2&0.62&0.13&0.43&-0.11 \\
		21.4&86.9&21.7&86.96&21.5&86.92&1.40&0.06&0.46&0.02 \\
		21.8&86.6&22.04&86.59&21.95&86.55&1.10&-0.01&0.68&-0.05 \\
		22.5&86&22.52&86.18&22.4&86.1&0.08&0.20&-0.44&0.11 \\
		22.8&85.8&23.06&85.75&22.75&85.7&1.14&-0.05&-0.21&-0.11 \\
		\hline
	\end{tabular}
\caption{Result of comparation with reference of TC, YAAS. A stands for error percentage between IMD and CIMSS; B stands for error percentage between IMD and Proposed Method}
\label{tab3}
\end{table}
\section{Conclusion} \label{sec5}
The purpose of this research study is to create a novel method for detecting the storm's eye utilising FIS and image processing. Based on expert knowledge acquired from the statistics of previous storm data, the suggested FIS pinpoints the location of the storm eye. The accuracy of the search is then improved through the use of image processing. The parameters are carefully analysed. Tauktae and Yaas, two TC, are used to test the methodology. The experimental findings lead to performance improvements of up to $98\%$. However, in the case of a high level of storm-scale, performance decline may happen and would be further enhanced in future research.

\section*{Declarations}
\subsection*{Conflict of Interest}
Conflict of interest is not declared by any of the authors.
\subsection*{Author Statement}
Research and manuscript preparation were done equally by each author.
\bibliography{sn-bibliography}


\begin{thebibliography}{46}
\ifx \bisbn   \undefined \def \bisbn  #1{ISBN #1}\fi
\ifx \binits  \undefined \def \binits#1{#1}\fi
\ifx \bauthor  \undefined \def \bauthor#1{#1}\fi
\ifx \batitle  \undefined \def \batitle#1{#1}\fi
\ifx \bjtitle  \undefined \def \bjtitle#1{#1}\fi
\ifx \bvolume  \undefined \def \bvolume#1{\textbf{#1}}\fi
\ifx \byear  \undefined \def \byear#1{#1}\fi
\ifx \bissue  \undefined \def \bissue#1{#1}\fi
\ifx \bfpage  \undefined \def \bfpage#1{#1}\fi
\ifx \blpage  \undefined \def \blpage #1{#1}\fi
\ifx \burl  \undefined \def \burl#1{\textsf{#1}}\fi
\ifx \doiurl  \undefined \def \doiurl#1{\url{https://doi.org/#1}}\fi
\ifx \betal  \undefined \def \betal{\textit{et al.}}\fi
\ifx \binstitute  \undefined \def \binstitute#1{#1}\fi
\ifx \binstitutionaled  \undefined \def \binstitutionaled#1{#1}\fi
\ifx \bctitle  \undefined \def \bctitle#1{#1}\fi
\ifx \beditor  \undefined \def \beditor#1{#1}\fi
\ifx \bpublisher  \undefined \def \bpublisher#1{#1}\fi
\ifx \bbtitle  \undefined \def \bbtitle#1{#1}\fi
\ifx \bedition  \undefined \def \bedition#1{#1}\fi
\ifx \bseriesno  \undefined \def \bseriesno#1{#1}\fi
\ifx \blocation  \undefined \def \blocation#1{#1}\fi
\ifx \bsertitle  \undefined \def \bsertitle#1{#1}\fi
\ifx \bsnm \undefined \def \bsnm#1{#1}\fi
\ifx \bsuffix \undefined \def \bsuffix#1{#1}\fi
\ifx \bparticle \undefined \def \bparticle#1{#1}\fi
\ifx \barticle \undefined \def \barticle#1{#1}\fi
\bibcommenthead
\ifx \bconfdate \undefined \def \bconfdate #1{#1}\fi
\ifx \botherref \undefined \def \botherref #1{#1}\fi
\ifx \url \undefined \def \url#1{\textsf{#1}}\fi
\ifx \bchapter \undefined \def \bchapter#1{#1}\fi
\ifx \bbook \undefined \def \bbook#1{#1}\fi
\ifx \bcomment \undefined \def \bcomment#1{#1}\fi
\ifx \oauthor \undefined \def \oauthor#1{#1}\fi
\ifx \citeauthoryear \undefined \def \citeauthoryear#1{#1}\fi
\ifx \endbibitem  \undefined \def \endbibitem {}\fi
\ifx \bconflocation  \undefined \def \bconflocation#1{#1}\fi
\ifx \arxivurl  \undefined \def \arxivurl#1{\textsf{#1}}\fi
\csname PreBibitemsHook\endcsname

\bibitem{bib3}
\begin{barticle}
\bauthor{\bsnm{Dvorak}, \binits{V.F.}}:
\batitle{Tropical cyclone intensity analysis and forecasting from satellite
  imagery}.
\bjtitle{Monthly Weather Review}
\bvolume{103}(\bissue{5}),
\bfpage{420}--\blpage{430}
(\byear{1975})
\end{barticle}
\endbibitem

\bibitem{bib4}
\begin{barticle}
\bauthor{\bsnm{Wong}, \binits{K.Y.}},
\bauthor{\bsnm{Yip}, \binits{C.L.}},
\bauthor{\bsnm{Li}, \binits{P.W.}}:
\batitle{Automatic identification of weather systems from numerical weather
  prediction data using genetic algorithm}.
\bjtitle{Expert systems with Applications}
\bvolume{35}(\bissue{1-2}),
\bfpage{542}--\blpage{555}
(\byear{2008})
\end{barticle}
\endbibitem

\bibitem{bib5}
\begin{bchapter}
\bauthor{\bsnm{Zou}, \binits{J.}},
\bauthor{\bsnm{Lin}, \binits{M.}},
\bauthor{\bsnm{Xie}, \binits{X.}},
\bauthor{\bsnm{Lang}, \binits{S.}},
\bauthor{\bsnm{Cui}, \binits{S.}}:
\bctitle{Automated typhoon identification from quikscat wind data}.
In: \bbtitle{2010 IEEE International Geoscience and Remote Sensing Symposium},
pp. \bfpage{4158}--\blpage{4161}
(\byear{2010}).
\bcomment{IEEE}
\end{bchapter}
\endbibitem

\bibitem{bib6}
\begin{bchapter}
\bauthor{\bsnm{Ho}, \binits{S.-S.}},
\bauthor{\bsnm{Talukder}, \binits{A.}}:
\bctitle{Automated cyclone identification from remote quikscat satellite data}.
In: \bbtitle{2008 IEEE Aerospace Conference},
pp. \bfpage{1}--\blpage{9}
(\byear{2008}).
\bcomment{IEEE}
\end{bchapter}
\endbibitem

\bibitem{bib7}
\begin{barticle}
\bauthor{\bsnm{Kumar}, \binits{P.}},
\bauthor{\bsnm{Kumar}, \binits{K.H.}},
\bauthor{\bsnm{Pal}, \binits{P.K.}}:
\batitle{Impact of oceansat-2 scatterometer winds and tmi observations on phet
  cyclone simulation}.
\bjtitle{IEEE transactions on geoscience and remote sensing}
\bvolume{51}(\bissue{6}),
\bfpage{3774}--\blpage{3779}
(\byear{2012})
\end{barticle}
\endbibitem

\bibitem{bib8}
\begin{barticle}
\bauthor{\bsnm{Zhong}, \binits{J.}},
\bauthor{\bsnm{Fei}, \binits{J.}},
\bauthor{\bsnm{Huang}, \binits{S.}},
\bauthor{\bsnm{Du}, \binits{H.}},
\bauthor{\bsnm{Zhang}, \binits{L.}}:
\batitle{An improved quikscat wind retrieval algorithm and eye locating for
  typhoon}.
\bjtitle{Acta Oceanologica Sinica}
\bvolume{31}(\bissue{1}),
\bfpage{41}--\blpage{50}
(\byear{2012})
\end{barticle}
\endbibitem

\bibitem{bib9}
\begin{barticle}
\bauthor{\bsnm{Jaiswal}, \binits{N.}},
\bauthor{\bsnm{Kishtawal}, \binits{C.M.}}:
\batitle{Prediction of tropical cyclogenesis using scatterometer data}.
\bjtitle{IEEE transactions on geoscience and remote sensing}
\bvolume{49}(\bissue{12}),
\bfpage{4904}--\blpage{4909}
(\byear{2011})
\end{barticle}
\endbibitem

\bibitem{bib10}
\begin{barticle}
\bauthor{\bsnm{Rodr{\'\i}guez}, \binits{P.G.}},
\bauthor{\bsnm{Polo}, \binits{M.-E.}},
\bauthor{\bsnm{Cuartero}, \binits{A.}},
\bauthor{\bsnm{Felic{\'\i}simo}, \binits{{\'A}.}},
\bauthor{\bsnm{Ruiz-Cuetos}, \binits{J.}}:
\batitle{Vecstatgraphs2d, a tool for the analysis of two-dimensional vector
  data: an example using quikscat ocean winds}.
\bjtitle{IEEE Geoscience and Remote Sensing Letters}
\bvolume{11}(\bissue{5}),
\bfpage{921}--\blpage{925}
(\byear{2013})
\end{barticle}
\endbibitem

\bibitem{bib11}
\begin{bchapter}
\bauthor{\bsnm{Lee}, \binits{R.S.}},
\bauthor{\bsnm{Liu}, \binits{J.N.}}:
\bctitle{Atmosphere-automatic track mining and objective satellite pattern
  hunting system using enhanced rbf and egdlm}.
In: \bbtitle{AAAI/IAAI},
pp. \bfpage{603}--\blpage{608}
(\byear{2000})
\end{bchapter}
\endbibitem

\bibitem{bib12}
\begin{barticle}
\bauthor{\bsnm{Lee}, \binits{R.}},
\bauthor{\bsnm{Liu}, \binits{J.}}:
\batitle{ijade weatherman: a weather forecasting system using intelligent
  multiagent-based fuzzy neuro network}.
\bjtitle{IEEE Transactions on Systems, Man, and Cybernetics, Part C
  (Applications and Reviews)}
\bvolume{34}(\bissue{3}),
\bfpage{369}--\blpage{377}
(\byear{2004})
\end{barticle}
\endbibitem

\bibitem{bib13}
\begin{barticle}
\bauthor{\bsnm{Terry}, \binits{J.P.}},
\bauthor{\bsnm{Feng}, \binits{C.-C.}}:
\batitle{On quantifying the sinuosity of typhoon tracks in the western north
  pacific basin}.
\bjtitle{Applied Geography}
\bvolume{30}(\bissue{4}),
\bfpage{678}--\blpage{686}
(\byear{2010})
\end{barticle}
\endbibitem

\bibitem{bib14}
\begin{barticle}
\bauthor{\bsnm{Bai}, \binits{Q.}},
\bauthor{\bsnm{Wei}, \binits{K.}}:
\batitle{Cloud system extraction in tropical cyclones by mountain-climbing}.
\bjtitle{Atmospheric Research}
\bvolume{101}(\bissue{3}),
\bfpage{611}--\blpage{620}
(\byear{2011})
\end{barticle}
\endbibitem

\bibitem{bib15}
\begin{barticle}
\bauthor{\bsnm{Bai}, \binits{Q.}},
\bauthor{\bsnm{Wei}, \binits{K.}},
\bauthor{\bsnm{Jing}, \binits{Z.}},
\bauthor{\bsnm{Li}, \binits{Y.}},
\bauthor{\bsnm{Tuo}, \binits{H.}},
\bauthor{\bsnm{Liu}, \binits{C.}}:
\batitle{Tropical cyclone spiral band extraction and center locating by binary
  ant colony optimization}.
\bjtitle{Science China Earth Sciences}
\bvolume{55}(\bissue{2}),
\bfpage{332}--\blpage{346}
(\byear{2012})
\end{barticle}
\endbibitem

\bibitem{bib16}
\begin{barticle}
\bauthor{\bsnm{Jaiswal}, \binits{N.}},
\bauthor{\bsnm{Kishtawal}, \binits{C.M.}}:
\batitle{Objective detection of center of tropical cyclone in remotely sensed
  infrared images}.
\bjtitle{IEEE Journal of Selected Topics in Applied Earth Observations and
  Remote Sensing}
\bvolume{6}(\bissue{2}),
\bfpage{1031}--\blpage{1035}
(\byear{2013})
\end{barticle}
\endbibitem

\bibitem{bib27}
\begin{barticle}
\bauthor{\bsnm{Lee}, \binits{R.S.}},
\bauthor{\bsnm{Liu}, \binits{J.N.}}:
\batitle{Tropical cyclone identification and tracking system using integrated
  neural oscillatory elastic graph matching and hybrid rbf network track mining
  techniques}.
\bjtitle{IEEE Transactions on Neural Networks}
\bvolume{11}(\bissue{3}),
\bfpage{680}--\blpage{689}
(\byear{2000})
\end{barticle}
\endbibitem

\bibitem{bib45}
\begin{bchapter}
\bauthor{\bsnm{Yan}, \binits{W.K.}},
\bauthor{\bsnm{Lap}, \binits{Y.C.}}:
\bctitle{An intelligent tropical cyclone eye fix system using motion field
  analysis}.
In: \bbtitle{17th IEEE International Conference on Tools with Artificial
  Intelligence (ICTAI'05)},
p. \bfpage{5}
(\byear{2005}).
\bcomment{IEEE}
\end{bchapter}
\endbibitem

\bibitem{bib17}
\begin{bchapter}
\bauthor{\bsnm{Warunsin}, \binits{K.}},
\bauthor{\bsnm{Chitsobhuk}, \binits{O.}}:
\bctitle{Automatic typhoon eye identification using quikscat data and spiral
  cloud image}.
In: \bbtitle{2014 14th International Conference on Control, Automation and
  Systems (ICCAS 2014)},
pp. \bfpage{212}--\blpage{216}
(\byear{2014}).
\bcomment{IEEE}
\end{bchapter}
\endbibitem

\bibitem{bib18}
\begin{bchapter}
\bauthor{\bsnm{Warunsin}, \binits{K.}},
\bauthor{\bsnm{Chitsobhuk}, \binits{O.}}:
\bctitle{Heuristic search on statistics of wind data and cloud images for
  automatic typhoon eye location}.
In: \bbtitle{2015 7th International Conference on Knowledge and Smart
  Technology (KST)},
pp. \bfpage{60}--\blpage{64}
(\byear{2015}).
\bcomment{IEEE}
\end{bchapter}
\endbibitem

\bibitem{bib19}
\begin{bchapter}
\bauthor{\bsnm{Warunsin}, \binits{K.}},
\bauthor{\bsnm{Chitsobhuk}, \binits{O.}}:
\bctitle{Heuristic search on statistics of wind data and cloud images for
  automatic typhoon eye location}.
In: \bbtitle{2015 7th International Conference on Knowledge and Smart
  Technology (KST)},
pp. \bfpage{60}--\blpage{64}
(\byear{2015}).
\bcomment{IEEE}
\end{bchapter}
\endbibitem

\bibitem{bib20}
\begin{barticle}
\bauthor{\bsnm{Lakshmanan}, \binits{V.}},
\bauthor{\bsnm{Rabin}, \binits{R.}},
\bauthor{\bsnm{DeBrunner}, \binits{V.}}:
\batitle{Multiscale storm identification and forecast}.
\bjtitle{Atmospheric research}
\bvolume{67},
\bfpage{367}--\blpage{380}
(\byear{2003})
\end{barticle}
\endbibitem

\bibitem{bib21}
\begin{bchapter}
\bauthor{\bsnm{Pao}, \binits{T.-L.}},
\bauthor{\bsnm{Yeh}, \binits{J.-H.}},
\bauthor{\bsnm{Liu}, \binits{M.-Y.}},
\bauthor{\bsnm{Hsu}, \binits{Y.-C.}}:
\bctitle{Locating the typhoon center from the ir satellite cloud images}.
In: \bbtitle{2006 IEEE International Conference on Systems, Man and
  Cybernetics},
vol. \bseriesno{1},
pp. \bfpage{484}--\blpage{488}
(\byear{2006}).
\bcomment{IEEE}
\end{bchapter}
\endbibitem

\bibitem{bib28}
\begin{barticle}
\bauthor{\bsnm{Harr}, \binits{P.A.}},
\bauthor{\bsnm{Elsberry}, \binits{R.L.}}:
\batitle{Large-scale circulation variability over the tropical western north
  pacific. part i: Spatial patterns and tropical cyclone characteristics}.
\bjtitle{Monthly Weather Review}
\bvolume{123}(\bissue{5}),
\bfpage{1225}--\blpage{1246}
(\byear{1995})
\end{barticle}
\endbibitem

\bibitem{bib29}
\begin{bchapter}
\bauthor{\bsnm{Fang}, \binits{Y.}},
\bauthor{\bsnm{Sugano}, \binits{K.}},
\bauthor{\bsnm{Oku}, \binits{K.}},
\bauthor{\bsnm{Kawagoe}, \binits{K.}}:
\bctitle{Applying a multi-dimensional time-series similarity method to
  typhoon-track prediction}.
In: \bbtitle{2015 IEEE 11th International Conference on e-Science},
pp. \bfpage{259}--\blpage{262}
(\byear{2015}).
\bcomment{IEEE}
\end{bchapter}
\endbibitem

\bibitem{bib30}
\begin{barticle}
\bauthor{\bsnm{Knaff}, \binits{J.A.}}:
\batitle{Revisiting the maximum intensity of recurving tropical cyclones}.
\bjtitle{International Journal of Climatology: A Journal of the Royal
  Meteorological Society}
\bvolume{29}(\bissue{6}),
\bfpage{827}--\blpage{837}
(\byear{2009})
\end{barticle}
\endbibitem

\bibitem{bib31}
\begin{bchapter}
\bauthor{\bsnm{Warunsin}, \binits{K.}},
\bauthor{\bsnm{Chitsobhuk}, \binits{O.}}:
\bctitle{Cyclone identification using fuzzy c mean clustering}.
In: \bbtitle{2013 13th International Symposium on Communications and
  Information Technologies (ISCIT)},
pp. \bfpage{369}--\blpage{373}
(\byear{2013}).
\bcomment{IEEE}
\end{bchapter}
\endbibitem

\bibitem{bib22}
\begin{barticle}
\bauthor{\bsnm{Pi{\~n}eros}, \binits{M.F.}},
\bauthor{\bsnm{Ritchie}, \binits{E.A.}},
\bauthor{\bsnm{Tyo}, \binits{J.S.}}:
\batitle{Objective measures of tropical cyclone structure and intensity change
  from remotely sensed infrared image data}.
\bjtitle{IEEE Transactions on Geoscience and Remote sensing}
\bvolume{46}(\bissue{11}),
\bfpage{3574}--\blpage{3580}
(\byear{2008})
\end{barticle}
\endbibitem

\bibitem{bib23}
\begin{barticle}
\bauthor{\bsnm{Pi{\~n}eros}, \binits{M.F.}},
\bauthor{\bsnm{Ritchie}, \binits{E.A.}},
\bauthor{\bsnm{Tyo}, \binits{J.S.}}:
\batitle{Estimating tropical cyclone intensity from infrared image data}.
\bjtitle{Weather and forecasting}
\bvolume{26}(\bissue{5}),
\bfpage{690}--\blpage{698}
(\byear{2011})
\end{barticle}
\endbibitem

\bibitem{bib24}
\begin{barticle}
\bauthor{\bsnm{Rodr{\'\i}guez-Herrera}, \binits{O.G.}},
\bauthor{\bsnm{Wood}, \binits{K.M.}},
\bauthor{\bsnm{Dolling}, \binits{K.P.}},
\bauthor{\bsnm{Black}, \binits{W.T.}},
\bauthor{\bsnm{Ritchie}, \binits{E.A.}},
\bauthor{\bsnm{Tyo}, \binits{J.S.}}:
\batitle{Automatic tracking of pregenesis tropical disturbances within the
  deviation angle variance system}.
\bjtitle{IEEE Geoscience and Remote Sensing Letters}
\bvolume{12}(\bissue{2}),
\bfpage{254}--\blpage{258}
(\byear{2014})
\end{barticle}
\endbibitem

\bibitem{bib25}
\begin{barticle}
\bauthor{\bsnm{Kossin}, \binits{J.P.}},
\bauthor{\bsnm{Velden}, \binits{C.S.}}:
\batitle{A pronounced bias in tropical cyclone minimum sea level pressure
  estimation based on the dvorak technique}.
\bjtitle{Monthly weather review}
\bvolume{132}(\bissue{1}),
\bfpage{165}--\blpage{173}
(\byear{2004})
\end{barticle}
\endbibitem

\bibitem{bib26}
\begin{bchapter}
\bauthor{\bsnm{Yan}, \binits{W.K.}},
\bauthor{\bsnm{Lap}, \binits{Y.C.}},
\bauthor{\bsnm{Wah}, \binits{L.P.}},
\bauthor{\bsnm{Wan}, \binits{T.W.}}:
\bctitle{Automatic template matching method for tropical cyclone eye fix}.
In: \bbtitle{Proceedings of the 17th International Conference on Pattern
  Recognition, 2004. ICPR 2004.},
vol. \bseriesno{3},
pp. \bfpage{650}--\blpage{653}
(\byear{2004}).
\bcomment{IEEE}
\end{bchapter}
\endbibitem

\bibitem{bib32}
\begin{barticle}
\bauthor{\bsnm{San~Pedro}, \binits{J.}},
\bauthor{\bsnm{Burstein}, \binits{F.}},
\bauthor{\bsnm{Sharp}, \binits{A.}}:
\batitle{A case-based fuzzy multicriteria decision support model for tropical
  cyclone forecasting}.
\bjtitle{European Journal of Operational Research}
\bvolume{160}(\bissue{2}),
\bfpage{308}--\blpage{324}
(\byear{2005})
\end{barticle}
\endbibitem

\bibitem{bib36}
\begin{barticle}
\bauthor{\bsnm{Badarinath}, \binits{K.}},
\bauthor{\bsnm{Kharol}, \binits{S.K.}},
\bauthor{\bsnm{Dileep}, \binits{P.}},
\bauthor{\bsnm{Prasad}, \binits{V.K.}}:
\batitle{Satellite observations on cyclone-induced upper ocean cooling and
  modulation of surface winds—a study on tropical ocean region}.
\bjtitle{IEEE Geoscience and Remote Sensing Letters}
\bvolume{6}(\bissue{3}),
\bfpage{481}--\blpage{485}
(\byear{2009})
\end{barticle}
\endbibitem

\bibitem{bib37}
\begin{barticle}
\bauthor{\bsnm{Hawkins}, \binits{J.D.}},
\bauthor{\bsnm{Turk}, \binits{F.J.}},
\bauthor{\bsnm{Lee}, \binits{T.F.}},
\bauthor{\bsnm{Richardson}, \binits{K.}}:
\batitle{Observations of tropical cyclones with the ssmis}.
\bjtitle{IEEE transactions on geoscience and remote sensing}
\bvolume{46}(\bissue{4}),
\bfpage{901}--\blpage{912}
(\byear{2008})
\end{barticle}
\endbibitem

\bibitem{bib38}
\begin{barticle}
\bauthor{\bsnm{Pi{\~n}eros}, \binits{M.F.}},
\bauthor{\bsnm{Ritchie}, \binits{E.A.}},
\bauthor{\bsnm{Tyo}, \binits{J.S.}}:
\batitle{Objective measures of tropical cyclone structure and intensity change
  from remotely sensed infrared image data}.
\bjtitle{IEEE Transactions on Geoscience and Remote sensing}
\bvolume{46}(\bissue{11}),
\bfpage{3574}--\blpage{3580}
(\byear{2008})
\end{barticle}
\endbibitem

\bibitem{bib39}
\begin{barticle}
\bauthor{\bsnm{Fritz}, \binits{S.}},
\bauthor{\bsnm{Hubert}, \binits{L.F.}},
\bauthor{\bsnm{Timchalk}, \binits{A.}}:
\batitle{Some inferences from satellite pictures of tropical disturbances}.
\bjtitle{Monthly Weather Review}
\bvolume{94}(\bissue{4}),
\bfpage{231}--\blpage{236}
(\byear{1966})
\end{barticle}
\endbibitem

\bibitem{bib40}
\begin{barticle}
\bauthor{\bsnm{Gentry}, \binits{R.C.}},
\bauthor{\bsnm{Fujita}, \binits{T.T.}},
\bauthor{\bsnm{Sheets}, \binits{R.C.}}:
\batitle{Aircraft, spacecraft, satellite and radar observations of hurricane
  gladys, 1968}.
\bjtitle{Journal of Applied Meteorology and Climatology}
\bvolume{9}(\bissue{6}),
\bfpage{837}--\blpage{850}
(\byear{1970})
\end{barticle}
\endbibitem

\bibitem{bib41}
\begin{barticle}
\bauthor{\bsnm{Rao}, \binits{B.}},
\bauthor{\bsnm{Kishtwal}, \binits{C.}},
\bauthor{\bsnm{Pal}, \binits{P.}},
\bauthor{\bsnm{Narayanan}, \binits{M.}}:
\batitle{Ers-1 surface wind observations over a cyclone system in the bay of
  bengal during november 1992}.
\bjtitle{Remote Sensing}
\bvolume{16}(\bissue{2}),
\bfpage{351}--\blpage{357}
(\byear{1995})
\end{barticle}
\endbibitem

\bibitem{bib42}
\begin{barticle}
\bauthor{\bsnm{Chaurasia}, \binits{S.}},
\bauthor{\bsnm{Kishtawal}, \binits{C.}},
\bauthor{\bsnm{Pal}, \binits{P.}}:
\batitle{An objective method of cyclone centre determination from geostationary
  satellite observations}.
\bjtitle{International Journal of Remote Sensing}
\bvolume{31}(\bissue{9}),
\bfpage{2429}--\blpage{2440}
(\byear{2010})
\end{barticle}
\endbibitem

\bibitem{bib43}
\begin{barticle}
\bauthor{\bsnm{Tuttle}, \binits{J.}},
\bauthor{\bsnm{Gall}, \binits{R.}}:
\batitle{A single-radar technique for estimating the winds in tropical
  cyclones}.
\bjtitle{Bulletin of the American Meteorological Society}
\bvolume{80}(\bissue{4}),
\bfpage{653}--\blpage{668}
(\byear{1999})
\end{barticle}
\endbibitem

\bibitem{bib44}
\begin{barticle}
\bauthor{\bsnm{Hasler}, \binits{A.}},
\bauthor{\bsnm{Palaniappan}, \binits{K.}},
\bauthor{\bsnm{Kambhammetu}, \binits{C.}},
\bauthor{\bsnm{Black}, \binits{P.}},
\bauthor{\bsnm{Uhlhorn}, \binits{E.}},
\bauthor{\bsnm{Chesters}, \binits{D.}}:
\batitle{High-resolution wind fields within the inner core and eye of a mature
  tropical cyclone from goes 1-min images}.
\bjtitle{Bulletin of the American Meteorological Society}
\bvolume{79}(\bissue{11}),
\bfpage{2483}--\blpage{2496}
(\byear{1998})
\end{barticle}
\endbibitem

\bibitem{bib35}
\begin{bchapter}
\bauthor{\bsnm{Liu}, \binits{J.N.-K.}},
\bauthor{\bsnm{Feng}, \binits{B.}},
\bauthor{\bsnm{Wang}, \binits{M.}},
\bauthor{\bsnm{Luo}, \binits{W.}}:
\bctitle{Tropical cyclone forecast using angle features and time warping}.
In: \bbtitle{The 2006 IEEE International Joint Conference on Neural Network
  Proceedings},
pp. \bfpage{4330}--\blpage{4337}
(\byear{2006}).
\bcomment{IEEE}
\end{bchapter}
\endbibitem

\bibitem{bib34}
\begin{barticle}
\bauthor{\bsnm{Yurchak}, \binits{B.S.}}:
\batitle{Description of cloud-rain bands in a tropical cyclone by a
  hyperbolic-logarithmic spiral}.
\bjtitle{Russian Meteorology and Hydrology}
\bvolume{32}(\bissue{1}),
\bfpage{8}--\blpage{18}
(\byear{2007})
\end{barticle}
\endbibitem

\bibitem{bib33}
\begin{barticle}
\bauthor{\bsnm{Wong}, \binits{K.Y.}},
\bauthor{\bsnm{Yip}, \binits{C.L.}},
\bauthor{\bsnm{Li}, \binits{P.W.}}:
\batitle{A novel algorithm for automatic tropical cyclone eye fix using doppler
  radar data}.
\bjtitle{Meteorological Applications: A journal of forecasting, practical
  applications, training techniques and modelling}
\bvolume{14}(\bissue{1}),
\bfpage{49}--\blpage{59}
(\byear{2007})
\end{barticle}
\endbibitem

\bibitem{bib1}
\begin{bchapter}
\bauthor{\bsnm{Riza}, \binits{L.S.}},
\bauthor{\bsnm{Bergmeir}, \binits{C.N.}},
\bauthor{\bsnm{Herrera~Triguero}, \binits{F.}},
\bauthor{\bsnm{Ben{\'\i}tez~S{\'a}nchez}, \binits{J.M.}}, \betal:
\bctitle{frbs: Fuzzy rule-based systems for classification and regression in
  r}.
(\byear{2015}).
\bcomment{American Statistical Association}
\end{bchapter}
\endbibitem

\bibitem{bib46}
\begin{botherref}
\oauthor{\bsnm{CIMSS}}:
CIMSS.
\url{http://tropic.ssec.wisc.edu/tropic.php}
\end{botherref}
\endbibitem

\bibitem{bib47}
\begin{botherref}
\oauthor{\bsnm{images}, \binits{S.}}:
Satellite images.
\url{http://satellite.imd.gov.in/insat.htm}
\end{botherref}
\endbibitem

\end{thebibliography}


\end{document}